\documentclass[10pt]{article}
\usepackage{amsmath}
\usepackage{amssymb}
\usepackage{amsthm}
\usepackage{a4wide}
\usepackage{amscd}
\usepackage{amsfonts}
\usepackage{graphicx}%
\usepackage{fancyhdr}
\usepackage{color}

\newcommand{\R}{\mathbb{R}}
\newcommand{\Om} {\Omega}
\newcommand{\pa} {\partial}
\newcommand{\be} {\begin{equation}}
\newcommand{\ee} {\end{equation}}
\newcommand{\bea} {\begin{eqnarray}}
\newcommand{\eea} {\end{eqnarray}}
\newcommand{\Bea} {\begin{eqnarray*}}
\newcommand{\Eea} {\end{eqnarray*}}

\newcommand{\al} {\alpha}

\newcommand{\de} {\delta}
\newcommand{\ga} {\gamma}

\newcommand{\De} {\Delta}
\newcommand{\la} {\lambda}

\newcommand{\rar}{\rightarrow}

\newcommand{\noi} {\noindent}

\newcommand{\va} {\varphi}

\newcommand{\e}{\epsilon}

\newcommand{\xs}{X^s_0(\Omega)}
\newcommand{\ds}{\displaystyle}
\newcommand{\ur}{\underline{u}}
\newcommand{\ts}{2_s^*}
\newcommand{\tss}{2_s^*-1}
\newcommand{\pp}{(\tilde{P}_\la^\mu)}

\theoremstyle{plain} \numberwithin{equation}{section}
\newtheorem{theorem}{Theorem}[section]

\newtheorem{lemma}[theorem]{Lemma}
\newtheorem{proposition}[theorem]{Proposition}
\theoremstyle{definition}
\newtheorem{definition}[theorem]{Definition}

\newtheorem{remark}[theorem]{Remark}

\begin{document}
 
\title{ A multiparameter semipositone fractional laplacian problem involving critical exponent}
 \author{R. Dhanya\thanks{School of Mathematics \& Computer Science, IIT Goa.
email: dhanya.tr@gmail.com
 }, Sweta Tiwari\thanks{ Department of Mathematics, IIT Guwahati. swetatiwari@iitg.ac.in}}

\maketitle
\begin{abstract}
In this paper we prove the existence of at least one positive solution for nonlocal semipositone problem of the type 
 $$
(P_\la^\mu)\left\{ \begin{array}{lll}
(-\De)^s u&=& \la(u^{q}-1)+\mu u^r \mbox{ in } \Om\\
u&>&0 \mbox{ in } \Om\\
u&\equiv &0 \mbox{ on }{\mathbb R^N\setminus\Om}.
\end{array}\right.
$$
when the positive parameters $\la$ and $\mu$ belongs to certain range. Here $\Om\subset \R^N$ is assumed to be a bounded open set  with smooth boundary, $s\in (0,1), N> 2s$ and $0<q<1<r\leq \frac{N+2s}{N- 2s}.$ 
The proof relies on the construction of a positive subsolution for $(P_\la^0)$ for  $\lambda>\la_0.$  
{
Now for each $\la>\la_0,$ for all small $0<\mu<\mu_{\la}$ we establish the existence of at least one positive solution of $(P_\la^\mu)$ using 
variational method. Also in the subcritical case, i.e., for $1<r<\frac{N+2s}{N-2s}$,
we show the existence of second
positive solution via mountain pass argument.
}
\end{abstract}, 
\section{Introduction}
Let $\Om\subset \R^N, $
  $ N > 2s$ be a bounded domain with boundary of class $C^{1,1}$ and we consider the following nonlocal partial differential equation
  with positive multiparameters $\la$ and $\mu $
  $$
(P_\la^\mu)\left\{ \begin{array}{lll}
(-\De)^s u&=& \la(u^{q}-1)+\mu u^r \mbox{ in } \Om\\
u&>&0 \mbox{ in } \Om\\
u&\equiv &0 \mbox{ on } {\mathbb R^N\setminus\Om}.
\end{array}\right.
$$
Here $0<q<1<r\leq \frac{N+2s}{N-2s}$ and $(-\De)^s$ denotes the standard fractional Laplace operator 
$$(-\De)^s u(x) =  2 C_{N,s}\mbox{ P.V.} \int_{\R^N}\frac{u(x)-u(y)}{|x-y|^{N+2s}} dy  $$
upto the renormalization factor $C_{N,s}$ given by
$\pi^{-\frac{N}{2}}2^{2s-1}\frac{\Gamma\left(\frac{N+2s}{2}\right)}{\Gamma(1-s)}$
where $\Gamma$ denotes the
Gamma function and  P.V. is understood in the sense of Cauchy's Principal Value. 
 Recently, great deal of attention has been given to the study of fractional and non-local operators of elliptic type,
 due to pure mathematical interests and also in the view of
concrete real-world applications. This type of operator arises in quite a natural
way in different contexts such as the thin obstacle problem,
optimization, finance, phase transitions,  anomalous diffusion etc. We refer 
\cite{Squass} and references there in for the problems involving non-local operator and their applications.

In this article we are interested to study the following non-local problem with semipositone nonlinearity.
  $$
\left\{ \begin{array}{lll}
(-\De)^s u&=& h(u)\mbox{ in } \Om\\
u&\equiv &0 \mbox{ on } \R^N\setminus \Om.
\end{array}\right.
$$
We say that the above non-local semilinear elliptic problem is semipositone in nature if 
$h(0)<0$ and $h$ is eventually positive. If $h$ is positive and monotone the above  problem is termed as  positone problem. 
Unlike positone problems where the positivity of non-negative solutions is guaranteed by the strong maximum principle, a semipositone problem can admit  
 non-negative solutions having zeros 
in the interior of $\Om$ even in the case of Laplacian, see \cite{Cas_Shi}.
Thus the mostly pursued existence theory by monotone iteration itself requires a positive subsolution which itself could be
challenging sometimes. Due to these reasons, in the celebrated works 
of Lions \cite{lions} and Berestycki\cite{BCN} it was mentioned that the semipositone problems are quite hard.  To the best of our knowledge the existence of positive solutions of semipositone problem 
for fractional Laplacian has not been studied so far. Here in this paper we address the  additional complexity which arises due to the 
critical exponent term $u^{\ts-1}.$

Our primary concern is to obtain a positive solution to the semipositone problem $(P_\la^0)$  
via monotone iteration(see section 3).
Motivated from the semipositone problem for Laplacian, the natural choice of 
subsolution of  $(P_\la^0)$ is a scalar multiple of $\phi_1^2$ where $\phi_1$ is the first eigenfunction of fractional Laplacian.  On the contrary to the standard Laplacian where the operator acts by pointwise differentiation, 
the fractional Laplace operators are defined via global integration and thus the proof of sub-solution is indeed much harder.  
In section 3 (Lemma \ref{subss}) of this paper we give a detailed calculation of this fact. Later in section 4, we use this particular function while 
defining a perturbed problem. It has come to our notice that in a recent work \cite{TGS}, authors obtain a weak solution for an infinite semipositone problem for a nonlocal problem using implicit function theorem .

In section 4 we are interested in finding positive solutions for the multiparameter problem $(P_\la^\mu)$ via variational method. 
In this regard we wish to mention that even in case of semilinear or quasilinear semipositone elliptic problems most of the existence
results were obtained using topological methods such as degree theory or bifurcation theory and sub-super solution methods (see \cite{DR} 
and the references therein).
A subcritical semipositone problem for Laplacian with one parameter is studied in \cite{CQY} using a non-smooth variational method.
In a more recent work (\cite{CQT}), the authors considers a more general class of nonlinearities with a
sign-changing weight via variational and continuity arguments. Positive solution of  semipositone super linear problem is established in the exterior of the ball for laplacian in  \cite{DQS} and for p-laplacian in \cite{QRI}. Recently in \cite{PSI}, multiplicity result for a 
multiparameter semipositone problem is obtained for a critical exponent problem via concentration compactness arguement.

In this paper we find solution(s) to the semi-positone problem via a variational approach similar to that of \cite{PS1}.
The idea is to define an appropriate cut-off functional  in such a way that the 
critical points on the newly defined functional itself becomes the weak solution of the original problem $(P_\la^\mu).$   
Here first critical point is obtained using direct methods of calculus of variations and the second solution, in case of sub-critical problem,
is obtained via mountain pass lemma.\\
Now we state our main results.
%
%
{First we consider the non-local problem with only semipositone term, $(P^0_\la)$ i.e  when $\mu=0.$
\begin{theorem} \label{T1}
There exists $\la_0\in(0, \infty)$ such that for all $\la>\la_0$ the problem $(P_\la^0)$ admits at least one positive solution and $(P_{\la_0}^0)$ admits a non-negative solution.
\end{theorem}
Next we consider the nonlinearities involving both semipositone and convex term and establish the following result: 

\begin{theorem}\label{T2}
For $1<r\leq\frac{N+2s}{N-2s}$ and for each $\la>\la_0$ ($\la_0$ as obtained in Theorem \ref{T1}),
there exists a $\mu_\la>0$ such that $(P_\la^\mu)$ 
admits at least one positive solution for all $\mu\in (0,\mu_\la)$. 
If  $1<r<\frac{N+2s}{N-2s}$,
there exists a second solution to the problem $(P_\la^\mu)$ for $\la>\la_0$ and  $\mu\in (0,{\mu}_\la)$.

\end{theorem}
}

\section
{Fractional Framework}

{In this section we shall discuss the functional setting to study the problem $(P_\la^\mu)$.
We refer \cite{BISCI} and \cite{hitch} for the details.}
\begin{definition}
The space $X^s(\Om)$ is the linear space of all Lebesgue measurable functions from 
$\R^N$ to $\R$ such that the restriction to $\Om$ of any function $u$ belongs to $L^2(\Om)$ and 
$$\displaystyle \int_{Q}\frac{|u(x)-u(y)|^2}{|x-y|^{N+2s}} dx\, dy<\infty \mbox{ where }  
Q=\R^N\times \R^N \setminus (\Omega^c \times \Omega^c ).$$ The space $X^s(\Om)$ is a Banach space equipped with the norm 
$$\|u\|_{X^s}=\|u\|_2+ \displaystyle \left( \int_{Q}\frac{|u(x)-u(y)|^2}{|x-y|^{N+2s}} dx \,dy\right)^{\frac{1}{2}}$$
{For addressing the Dirichlet's boundary condition in this paper, we will now consider the following linear subspace of $X^s(\Om)$.}
$$X^s_0(\Om)=\{u \in X^s(\Om) : u =0\,\, a.e \mbox{ in } \Om^c \}.$$
\end{definition}
\noindent Note that $C_c^\infty(\Om)$ is dense in $X_0^s(\Om)$ and $X^s_0(\Om)$ is a Hilbert space equipped with the inner product 
$$(u,v) = \int _{\Om} u v +\int_{Q}\frac{(u(x)-u(y))(v(x)-v(y))}{|x-y|^{N+2s}}dx \, dy.$$
It is also known that $X^s_0(\Om)=\{u\in H^s(\R^N): u =0 \, a.e \, \mbox{ in } \Om^c\}$ where $H^s(\R^N)$ denote 
the fractional order Sob\"{o}lev space in $\R^N$ defined by the Fourier transform. We have the following continuous embedding 
$H^s(\R^N) \hookrightarrow L^{2_s^*}(\R^N)$ where $2_s^*$ is the critical Sobolev exponent $\frac{2N}{N-2s}.$ In fact there exists a constant $C>0$ such that  
$$\|u\|_{L^{2_s^*}(\R^N)}\leq C\int_{\Om} \frac{|u(x)-u(y)|^2}{|x-y|^{N+2s}} dx\,dy \;\;\; \mbox{ holds for all } u\in H^s(\R^N)$$
and as a consequence of which  we conclude that $X^s_0(\Om) \hookrightarrow L^{2_s^*}(\R^N)$ and since $\Om$ is a bounded open set in $\R^N$
$$u \mapsto \left(\int_{Q}\frac{|u(x)-u(y)|^2}{|x-y|^{N+2s}} dx\,dy\right)^{\frac{1}{2}} $$
defines the equivalent norm in $X_0^s(\Om).$ 
Throughout this paper by the term $\|u\|^2$ we always mean that $\ds \int_{Q}\frac{|u(x)-u(y)|^2}{|x-y|^{N+2s}} dx\,dy$. 
\begin{definition}
We say $u\in \xs$ is a weak solution for $(P^\mu_\la)$ if the following identity  holds for all $\va\in C_c^\infty(\Om)$ 
$$C_{N,s}\int_Q  \frac{(u(x)-u(y))(\va(x)-\va(y))}{|x-y|^{N+2s}}dx \,dy = \la \int_{\Om} (u^q-1)\va +\mu \int_{\Om} u^r \va. $$
\end{definition}

\section{Sub-supersolution approach for  semipositone problem}
In this section we give the proof of Theorem \ref{T1}.
We aim to prove the existence of a positive solution of the following problem, $({P}_\la^0)$, via monotone iteration technique.
$$
(P_\la^0 )\left\{ \begin{array}{rll}
(-\De)^s u&=& \la( u^{q}-1) \mbox{ in } \Om,\\
u&>&0 \mbox{ in } \Om,\\
u&\equiv &0 \mbox{ in }\R^N \setminus \Om.
\end{array}\right.
$$
 The crucial step 
towards this is to find a suitable positive subsolution . So far we are not aware of such a treatment for a fractional Laplace
problem with a semipositone nonlinearity. In the next lemma we construct a weak subsolution for $(P_\la^0)$ and then use weak comparison principle to obtain a monotone sequence of functions which converges to a weak solution to $(P_\la^0).$

Let $\la_1$ denotes the least eigenvalue of $(-\De)^s$ in $\Om$ and $\phi_1$ be the corresponding eigenfunction, i.e 
$$
\left\{ \begin{array}{rll}
(-\De)^s \phi_1&=& \la_1 \phi_1 \mbox{ in } \Om,\\
\phi_1&\equiv &0 \mbox{ in }\R^N \setminus \Om.
\end{array}\right.
$$
\begin{lemma}\label{subss}
The function $\underline{u}= \la^{\alpha_1}\phi_1^2$ where $\alpha_1\in (1,\frac{1}{1-q})$ is a {weak}
subsolution of $(P_\la^0)$ for $\lambda$ large. 
\end{lemma}
\noi Proof: We know that $\phi_1$ is bounded and strictly positive inside $\Om$ (for details see  chapter 3 of \cite{BISCI}).
Since $\phi_1\in L^\infty(\R^N)\cap X_0^s(\Om),$ we observe that the function $\phi_1^2$ belongs to $ X_0^s(\Om).$  Next we compute $(-\De)^s\phi_1^2(x)$ for $x\in \Omega$ using the definition of the fractional Laplacian.  For each $\e>0$ and $x\in\Om$ we write 

$$ \int_{B_\e (x)^c}\frac{\phi_1^2(x)-\phi_1^2(y)}{|x-y|^{N+2s} }\, dy = 2 \phi_1(x) \int_{B_\e(x)^c} \frac{\phi_1(x)-\phi_1(y)}{|x-y|^{N+2s}} -\int_{B_\e(x)^c} \frac{(\phi_1(x)-\phi_1(y))^2}{|x-y|^{N+2s}}. $$
Clearly for each $x\in \Om,\,$  $\displaystyle \lim_{\e\rar 0} \int_{B_\e(x)^c} \frac{\phi_1(x)-\phi_1(y)}{|x-y|^{N+2s}}$ exists and is equal to $\displaystyle\frac{\la_1 \phi_1(x)}{2 C_{N,s}}.$ Letting $\e\rar 0$ in the previous expression we have
$$\begin{array}{lll}

\ds \mbox{P.V.} \int_{\R^N}\frac{\phi_1^2(x)-\phi_1^2(y)}{|x-y|^{N+2s} }\, dy& = \ds 2 \phi_1(x) .\, \, \mbox{P.V.} \int_{\R^N} \frac{\phi_1(x)-\phi_1(y)}{|x-y|^{N+2s}} - \mbox{P.V.} \int_{\R^N} \frac{(\phi_1(x)-\phi_1(y))^2}{|x-y|^{N+2s}}\\[3mm]
&= \ds \frac{1}{C_{N,s}}   \phi_1 (-\De)^s \phi_1-\mbox{ P.V.} \int_{\R^N} \frac{(\phi_1(x)-\phi_1(y))^2}{|x-y|^{N+2s}} dy.
\end{array}
$$
Since the integrand is positive, P.V $ \ds  \int_{\R^N} \frac{(\phi_1(x)-\phi_1(y))^2}{|x-y|^{N+2s}} dy = \int_{\mathbb R^N}\frac{(\phi_1(x)-\phi_1(y))^2}{|x-y|^{N+2s}}dy.$ Now if we write 
\begin{equation}\label{Ispositive} h(x)=\int_{\mathbb R^N}\frac{(\phi_1(x)-\phi_1(y))^2}{|x-y|^{N+2s}}dy  \mbox{ for all } x\in \Omega  
\end{equation}
then clearly we have 
\begin{equation}\label{eqn3.2}
(-\De)^s \phi_1^2(x) = 2 \la_1 \phi_1^2    - 2 C_{N,s}h(x) \mbox{ for all } x\in \Om.
\end{equation} Next we claim that $0<h(x)<\infty$ for all $x\in \Om.$ 
 To prove the claim we use the interior regularity result for  $\phi_1,$ i.e the first eigenfunction $\phi_1$ belongs to $C^{2 s+\alpha}$ inside $\Om$ (Theorem 6.1 of see \cite{oton_survey}). Now if $0<s<\frac{1}{2},$ we first choose $\alpha$ small enough so that $2s+\alpha<1.$   Now, let $x\in \Om$ and fix an $r<<1$ so that $B_r(x)$ is compactly contained in $\Om.$ Now by interior regularity result $|\phi_1(x)-\phi_1(y)|\leq C|x-y|^{2s
 +\alpha}$ for all $y\in B_r(x).$ Thus
 \begin{equation} 
\int_{\R^N}\frac{(\phi_1(x)-\phi_1(y))^2}{|x-y|^{N+2s}}dy \leq 2\, C \,\|\phi_1\|_{\infty} \int_{|x-y|<r}  |x-y|^{-N+ \alpha} \,dy+ 4\, \|\phi_1 \|_{\infty}^2 \int_{|x-y|>r} \frac{1}{|x-y|^{N+2s}}\,dy 
\end{equation}
Integrating in the polar co-ordinates we find that $0<h(x)<\infty$ for all $x\in \Om$ and $s\in (0,\frac{1}{2}).$ Now if $s\in [\frac{1}{2},1),$ the first eigen function is clearly Lipschitz continuous inside $\Om$ and hence 
\begin{equation} 
\int_{\R^N}\frac{(\phi_1(x)-\phi_1(y))^2}{|x-y|^{N+2s}}dy \,\leq  C  \int_{|x-y|<r}  |x-y|^{-N-2s +2} \,dy+ 4\|\phi_1 \|_{\infty}^2 \int_{|x-y|>r} \frac{1}{|x-y|^{N+2s}}\,dy < \infty.
\end{equation}
Next observe that for all $x\in \Om,$ $h(x)>\delta $ for some $\de>0.$ Indeed, if $\{x_n\}$ be a sequence in $ \Om$ such that  
$h(x_n)\rar 0$  and  $x_{n}\rar x_0$ for some $x_0\in \overline{\Om}$, then by  Fat\"{o}u's Lemma, $h(x_0)=0$ which contradicts the fact that $h(x)>0$
in $\overline{\Om}.$ Also note that the function $h(x)$ extended upto $\pa\Om$ may possibly assume the value $+\infty$ as $x\rar \pa\Om.$ 
Now we shall prove that $\underbar{u}=\la^{\alpha_1} \phi_1^2$ is a subsolution of $(P_\la^0).$   
Since $\alpha_1>1$ and $h(x)>\de_0$ in $\overline{\Omega}$  for some $\de_0>0$, we have for large values of $\lambda$ 
\begin{equation}
\label{subsol33}
2 C_{N,s} \la^{\al_1}h(x)>\la \text{ for all } x\in {\Omega}.
 \end{equation}
Also as $\phi_1\in L^\infty(\R^N)$ and $\alpha_1<\frac{1}{1-q}$, again for large  $\la$, we have
\begin{equation}\label{subsol44}
2\lambda^{\alpha_1}\la_1\phi_1^2\leq \la^{\alpha_1q+1}\phi_1^{2q}\text{ for } x\in\Om.
\end{equation} 
Thus from \eqref{eqn3.2} and \eqref{subsol33}-\eqref{subsol44},
we find that for large values of $\la$ and for $x\in\Om$,
\begin{equation}\label{subsol55}
(-\De)^s \underline{u}(x) = 2\lambda^{\alpha_1}\la_1\phi_1^2(x)-2 C_{N,s} \lambda^{\alpha_1}h(x)<\la (\underline{u}^q-1).
\end{equation}
Next it remains to prove that $\underline{u}$ is a weak subsolution of $(P_\la^0)$, i.e.,
$$C_{N,s}\int_Q  \frac{(\underline{u}(x)-\underline{u}(y))(\va(x)-\va(y))}{|x-y|^{N+2s}}dx \,dy \leq \la \int_{\Om} (\ur^q-1)\va \mbox{ holds for all } \va\in C_c^\infty(\Om) , \,  \va\geq 0.$$ 
We verify the claim in two steps. \\
\underline{Step I:} First we show that
$$\int_Q  \frac{(\underline{u}(x)-\underline{u}(y))(\va(x)-\va(y))}{|x-y|^{N+2s}}dx \,dy= \int_{\R^N\times \R^N} \va(x) \left(\frac{2\underline{u}(x)-\underline{u}(x+\xi)-\underline{u}(x-\xi)}{|\xi |^{N+2s}}\right) dx \,d\xi.$$
 To verify the step I we use the same idea of the proof of integration by parts formula for fractional Laplacian as in Lemma 1.26 of \cite{BISCI}. We can define the sets $\mathcal{D}_0, \mathcal{D}_\e, \mathcal{D}_\e^{\pm}$ as in Lemma 1.26. It is also easy to verify that 
 \begin{equation}\label{step1eq}
 \int_{\mathcal{D_\e}} \frac{(\underline{u}(x)-\underline{u}(y))(\va(x)-\va(y))}{|x-y|^{N+2s}}dx \,dy= \int_{\mathcal{D_\e^{+}}\cup \mathcal{D_\e^{-}}} \va(x) \left(\frac{2\underline{u}(x)-\underline{u}(x+\xi)-\underline{u}(x-\xi)}{|\xi |^{N+2s}}\right) dx d\xi.
 \end{equation}
 Since $\underline{u}$ and $\va$ belongs to $X_0^s(\Om)$, letting $\e\rar 0$ in the left hand side of the above equation, we obtain 
 \begin{equation}
 \int_{\mathcal{D_\e}} \frac{(\underline{u}(x)-\underline{u}(y))(\va(x)-\va(y))}{|x-y|^{N+2s}}dx \,dy=\int_{Q} \frac{(\underline{u}(x)-\underline{u}(y))(\va(x)-\va(y))}{|x-y|^{N+2s}}dx \,dy.
 \end{equation}
 Define $U(x,\xi)= \ds \left(\frac{2\underline{u}(x)-\underline{u}(x+\xi)-\underline{u}(x-\xi)}{|\xi |^{N+2s}}\right).$ 
 Though it is not clear if $U(x,\xi) \in L^1(\R^N\times \R^N)$, we have $\va(x) U(x,\xi)\in L^{1}(\R^N\times \R^N)$. Indeed, as
 $\va$ has support inside $\Om_0$ which is compactly contained in $\Om$, we have \\
 $\ds \int_{\mathcal{D_\e^{+}}\cup \mathcal{D_\e^{-}}} \va(x) \left(\frac{2\underline{u}(x)-\underline{u}(x+\xi)-\underline{u}(x-\xi)}{|\xi |^{N+2s}}\right)$
 $$ \;\;\;\;\;\;\;\;\;\;=\int_{\{ x\in \Om_0\} \cap \mathcal{D_\e^{+}}\cup \mathcal{D_\e^{-}}} \va(x) \left(\frac{2\underline{u}(x)-\underline{u}(x+\xi)-\underline{u}(x-\xi)}{|\xi |^{N+2s}}\right) .$$
 Let $ 2 r:=dist(\Om_0,\pa\Om)$ so that for each $x\in \Om_0$ the ball $B_r(x)$ is compactly contained inside $\Om.$ 
 Let $\Om_1=\cup_{x\in \Om_0} B_r(x),$ then once again we use the interior regularity result to conclude that 
 $\phi_1\in C^2(\overline{\Om_1}).$ Thus by Taylor's theorem we have
\begin{equation}
 \left|\frac{2\underline{u}(x)-\underline{u}(x+\xi)-\underline{u}(x-\xi)}{|\xi |^{N+2s}}\right| \leq \|D^2\underline{u}\|_{L^{\infty}{(\Om_1})} |\xi|^{2-2s-N} \,\,\,\,  \mbox{ for }\xi\in B_r(0)
\end{equation}
 and 
 \begin{equation}
 \left|\frac{2\underline{u}(x)-\underline{u}(x+\xi)-\underline{u}(x-\xi)}{|\xi |^{N+2s}}\right| \leq 4 \|\underline{u}\|_{\infty} |\xi|^{-2s-N} \,\,\,\,  \mbox{ for } |\xi|>r.
 \end{equation} 
 Now we evaluate $ \ds \int_{x\in \Om_0}\int_{\xi\in \R^N}  \va(x) \left|\frac{2\underline{u}(x)-\underline{u}(x+\xi)-\underline{u}(x-\xi)}{|\xi |^{N+2s}}\right|$
$$
\begin{array}{lll}
&=& \ds \int_{x\in \Om_0}\int_{\xi\in B_r(0)}  \va(x) \left|\frac{2\underline{u}(x)-\underline{u}(x+\xi)-\underline{u}(x-\xi)}{|\xi |^{N+2s}}\right|\\[4mm]
&& \ds  + \int_{x\in \Om_0}\int_{\xi\in (B_r(0))^c}  \va(x) \left|\frac{2\underline{u}(x)-\underline{u}(x+\xi)-\underline{u}(x-\xi)}{|\xi |^{N+2s}}\right|\\[4mm]
&\leq& \ds \int_{x\in \Om_0} \int_{\xi\in B_r(0)}  \|D^2\underline{u}\|_{L^{\infty}{(\Om_1})} |\xi|^{2-2s-N} d\xi dx + \\[4mm]
&&+\ds \int_{x\in \Om_0}\ \int_{\xi\in( B_r(0))^c} 4 \|\underline{u}\|_{\infty} |\xi|^{-2s-N} d\xi dx.\\
&<& \infty.
\end{array} 
 $$
 This implies that $\va(x) U(x,\xi)\in L^{1}(\R^N\times \R^N)$ and hence we can pass through the limit in the right hand side of \eqref{step1eq} and 
 thus the Step I is verified. \\
\underline{Step II:}  We claim that
 $$\int_{\R^N\times \R^N} \va(x) \left(\frac{2\underline{u}(x)-\underline{u}(x+\xi)-\underline{u}(x-\xi)}{|\xi |^{N+2s}}\right) dx \,d\xi= \frac{1}{C_{N,s}} \int_{\R^N} \va(x) (-\Delta)^s \underline{u}(x) dx. $$
First note that the $\lim_{\e\rar 0}$
$\ds \int_{\R^N\setminus B_\e(y)} \frac{\underline{u}(x)- \underline{u}(y)}{|x-y|^{N+2s} } dy$ exists for all $x\in \Om,$ and hence we can write 

\begin{eqnarray}
 \ds 2\, \mbox{P.V}  \int_{\R^N}\frac{\underline{u}(x+z)- \underline{u}(x)}{|z|^{N+2s}} dz &=& 
 \ds \mbox{P.V}  \int_{\R^N}\frac{\underline{u}(x+z)- \underline{u}(x)}{|z|^{N+2s} } dz \nonumber \\[3mm]
 & & +  \mbox{P.V} \ds \int_{\R^N}\frac{\underline{u}(x-z)- \underline{u}(x)}{|z|^{N+2s} } dz\nonumber  \\[3mm]
 &=&  \mbox{P.V}  \ds \int_{\R^N}\frac{\underline{u}(x+y)+ \underline{u}(x-y)- 2 \underline{u}(x)}{|y|^{N+2s} } dy.\label{eqn312}
\end{eqnarray} 
 Now from Step I we can apply Fubini's theorem and then use (\ref{eqn312}) to obtain\\[3mm]
 $ \ds \int_{\R^N\times \R^N} \va(x) \left(\frac{2\underline{u}(x)-\underline{u}(x+y)-\underline{u}(x-y)}{|y |^{N+2s}}\right) dx \,dy  $ 
 $$
 \begin{array}{lll} 
 &=& \ds \int_{ \Om_0} \va(x)\int_{ \R^N} \left(\frac{2\underline{u}(x)-\underline{u}(x+y)-\underline{u}(x-y)}{|y |^{N+2s}}\right) dy \,dx .\\[4mm]
 &=& \ds \int_{x\in \Om_0} \va(x) (-2)\, \mbox{P.V}  \int_{\R^N}\left( \frac{\underline{u}(x+y)- \underline{u}(x)}{|y|^{N+2s}} \right) dy\, dx\\[3mm]
 &=& \ds \frac{1}{C_{N,s}} \int_{\R^N} \va(x) (-\Delta)^s \underline{u}(x) dx.
 
\end{array} 
 $$
which completes the proof of Step II and hence the Lemma is proved.\hfill\qed

 \begin{theorem}
There exists $\la_0\in(0, \infty)$ such that for all $\la>\la_0$ the problem $(P_\la^0)$ admits at least one positive solution and $(P_{\la_0}^0)$ admits a non-negative solution.
\end{theorem}
\noi Proof: Consider the unique weak solution $z_\la \in X^s_0(\Omega)$ of the non-local problem
$${(\bar{P}_{\la})}
\left\{ \begin{array}{rll}
(-\De)^s z_\la&=& \la  z_{\la}^{q} \mbox{ in } \Om,\\
z_\la&>&0 \mbox{ in } \Om,\\
z_\la &\equiv &0 \mbox{ in }\R^N \setminus \Om.
\end{array}\right.
$$
The existence of $z_\la$ can be easily proved via standard minimization of the associated functional
$E_1^\la(u) = \displaystyle\frac{C_{N,s}}{2}  \int_{Q} \frac{|u(x)-u(y)|^2}{|x-y|^{N+2s}} -\frac{\la}{q+1}\int_{\Om} |u|^{q+1}$ in $\xs.$
Now by {Theorem 3.2 in \cite{FP}}, $z_\la$ is in $L^\infty(\Om)$. Using the regularity result in \cite{Rosoton1},
$z_\la\in C^s(\R^N)$ and $u/d(x)^s\in C^{\alpha}(\overline{\Om})$ for some $\alpha>0.$
Clearly we have $z_\la= \la^{\frac{1}{1-q}}z_1$  and by Hopf' s Lemma  (see Lemma 7.3 of \cite{oton_survey})
$z_\la(x) \geq C\, \la^{\frac{1}{1-q}} d(x)^s$ in $\Om.$ Clearly $z_\la$ is a super solution of {$(P_\la^0)$} for all $\la>0.$
From Lemma \ref{subss}, $\ur$ is a subsolution of
$({P}_\la^0)$ for $\la>>1.$ Now that both $\phi_1$ and $z_\la$ behaves like $(d(x))^s$ inside $\Omega,$
if required we can choose $\la$ still larger so that $\underline{u}\leq z_\la.$
Then by the standard monotone iteration method for
the fractional Laplacian, we obtain a solution $u_{\la}$ of $(P_\la^0)$ lying in between the ordered pair
$[\underline{u},  z_\la].$
 Define the set $\Lambda:=\{\la>0:(P_\la^0) \text{ admits a weak solution}\}$ and let $\la_0=\inf\Lambda$.
 Fix a $\la>\la_0$ and by the definition of $\la_0$ there exists a $\la'\in (\la_0,\la)$ such that $(P_{\la'}^0)$ admits a weak solution $u'.$
 { Now we can proceed as in the proof of Theorem 2.1 of \cite{DR} and obtain a solution $u_\la$ which lies in between the ordered interval $[u',z_\la].$ Now combining both the steps we assert
that the problem $(P_\la^0)$ admits one positive solution
all $\la>\la_0$.}

 Next we show the existence of a non-negative weak solution for the problem $(P_{\la_0}^0)$. Consider a minimizing sequence $\{\la_n\}\subset\Lambda$ converging to $\la_0.$  Then as $\la_n\downarrow\la_0$, we get  
$\tilde \lambda\in\Lambda$ such that
$\la_n\leq \tilde\lambda$ for all $n\in\mathbb N$. Let $\{u_{\la_n}\}$, $\{z_{\la_n}\}$ and $z_{\tilde\la}$ be the solutions of the problems $(P_{\la_n}^0)$,$(\bar{P}_{\la_n}^0)$ and $(\bar{P}_{\tilde\la}^0)$ respectively.
By the construction of the solution, we have $0<u_{\la_n}\leq z_{\la_n}\leq z_{\tilde\la}$ for all $n\in\mathbb N$ in $\Om$.
Also as \begin{align*}
         \int_Q  \frac{(u_{\la_n}(x)-u_{\la_n}(y))^2}{|x-y|^{N+2s}}dx \,dy=\la_n \int_{\Om} (u_{\la_n}^q-1)u_{\la_n}
         \leq\la_n \int_{\Om} (z_{\la_n}^q+1)z_{\la_n}
        \leq\tilde\la\int_{\Om} (z_{\tilde\la}^q+1)z_{\tilde\la}
        \end{align*}
it implies that the sequence $\{u_{\la_n}\}$ is  bounded in $X^s_0(\Om)$ and hence converges weakly in $X_0^s(\Om)$ to $u_{\la_0}$ (say).
Then, passing through the limit in the weak sense, one can easily check that $u_{\la_0}$ is a non-negative weak solution of $(P_{\la_0}^0)$.
This completes the Theorem 3.2. \hfill\qed

{\begin{remark}
Here we note that the sup-super solution method can also be employed for finding a solution of the problem $(P_\la^\mu)$ for $\mu>0$ and small. 
Indeed from Lemma \ref{subss},  $\underline{u}= \la^{\alpha_1}\phi_1^2$ where $\alpha_1\in (1,\frac{1}{1-q})$, is a weak
subsolution of $(P_\la^\mu)$ for $\lambda$ large and all $\mu>0$. Also for $\al_2>\frac{1}{1-q}$, one can check that $\overline{u}=\la^{\al_2}\psi$
is a super-solution for $(P_\la^\mu)$ for $\la$ large and $\mu>0$ small enough, where $\psi$ solves
$$
\left\{ \begin{array}{rll}
(-\De)^s \psi&=& 1 \mbox{ in } \Om,\\
\psi&>&0 \mbox{ in } \Om,\\
\psi &\equiv &0 \mbox{ in }\R^N \setminus \Om.
\end{array}\right.
$$
\end{remark}}

\section{Variational approach towards concave-convex problem}
In this section we study the combined effect of the semipositone and convex terms on the existence and multiplicity of the solutions for the problem $(P_\la^\mu)$ for $\mu>0$.
Now for $0<q<1$ and $1< r \leq  2_s^*-1,$ define the cut-off problem 
 $$
(\tilde{P}_{\la}^\mu)\left\{ \begin{array}{lll}
(-\De)^s u&=& \la g(x,u)+\mu f(x,u) \mbox{ in } \Om\\
u&>&0 \mbox{ in } \Om\\
u&\equiv &0 \mbox{ on } \R^N\setminus \Omega.
\end{array}\right.
$$
where 
$$
f(x,t)=\left\{ \begin{array}{lll}
\underline{u}^r(x) & \mbox{ if }  t\leq \underline{u}(x)\\
t^r& \mbox{ if } \underline{u}(x)<t.\\
\end{array}\right.
$$
and 
$$
g(x,t)=\left\{ \begin{array}{lll}
\underline{u}^q(x) -1& \mbox{ if } t\leq  \underline{u}(x)\\
t^q-1& \mbox{ if } \underline{u}(x)<t.\\
\end{array}\right.
$$
and $\underline{u}$ is as defined in section 3. Also define $F(x,t)=\int_0^t f(x,s)ds$ and $G(x,t)=\int_0^t g(x,s)ds.$  Then explicitly we can write
$$
F(x,t)=\left\{ \begin{array}{lll}
\ds \underline{u}^r t & \mbox{ if } t\leq  \underline{u}(x)\\[4mm]
\ds \frac{t^{r+1}}{r+1} +\frac{r}{r+1} \underline{u}^{r+1}& \mbox{ if } \underline{u}(x)<t.\\
\end{array}\right.
$$
and 
$$
G(x,t)=\left\{ \begin{array}{lll}
\ds (\ur^q-1) t & \mbox{ if } t\leq  \underline{u}(x)\\[4mm]
\ds \frac{t^{q+1}}{q+1}-t+\frac{q}{q+1} \ur^{q+1}.& \mbox{ if } \underline{u}(x)<t.\\
\end{array}\right.
$$

A direct calculation yields that there exists positive constants $c, c'$ such that 
\begin{equation}\label{estFG}
|F(x,t)|\leq c+c'|t|+ \frac{|t|^{r+1}}{r+1} \mbox{ and  }|G(x,t)|\leq c+c' |t|+ \frac{|t|^{q+1}}{q+1}  
\end{equation}

Let $I_\mu$ denote the energy functional associated to cut-off problem given by
\begin{equation}\label{Id}
I_\mu(u) =\frac{1}{2}\|u\|^2-\mu\int_{\Om} F(x,u)-\la \int_{\Om} G(x,u) \mbox{ for } u\in \xs.
\end{equation}
{Positivity of solution(s) of $(P_\la^\mu)$ is obtained using the following important remark. }
\begin{remark}\label{rem4.4}
 $I_\mu$ is a $C^1$ functional on $\xs$ and the critical points of $I_\mu$ are weak solutions of $\pp.$ By weak comparison principle for fractional laplacian any such weak solution is bounded below by $\ur$ and hence positive weak solution of $(P_\la^\mu)$ itself.
\end{remark}

\begin{lemma}\label{PSC}
Let $I_\mu$ be defined as in $\eqref{Id}$ and $1< r<\ts-1.$ Then $I_{\mu}$ satisfies $ (PS)_c$ condition for every $c\in \R.$
\end{lemma}
  \noi Proof: 
 Let $\{u_n\}$ be a Palais Smale sequnce at a level $c,$ i.e $I_{\mu}(u_n)\rar c$ and $I_{\mu}'(u_n)\rar 0.$ First we show 
 that $\{u_n\}$ is bounded in 
 $X_0^s(\Om)$. For this set $A_n=\{x: u_n >\ur\}$ and $B_n=\{x: u_n\leq \ur\}.$ Then, 
\begin{equation}\label{PSbdd1}
\begin{array}{lll}
 c+o_n(1)&=&\ds I_\mu(u_n)- \frac{< I_\mu'(u_n),u_n >}{r+1}\\[3mm]
 &=&\ds\left( \frac{1}{2}-\frac{1}{r+1}\right)\|u_n\|^2-\mu\left(\int_\Om F(x,u_n)-\frac{f(x,u_n)u_n}{r+1}\right)
 -\la\left(\int_\Om G(x,u_n)-\frac{g(x,u_n)u_n}{r+1}\right)\\[3mm]

 \end{array}
 \end{equation}
 Now we bound the second and third integrals in terms of $\|u_n\|$ as below: 
 \begin{equation}\label{PSbdd2} 
 \begin{array}{lll}
 \ds\left|\int_\Om F(x,u_n)-\frac{f(x,u_n)u_n}{r+1}\right|&=&\left|\ds\left(\int_{A_n} F(x,u_n)-\frac{f(x,u_n)u_n}{r+1}\right)+
 \ds\left(\int_{B_n} F(x,u_n)-\frac{f(x,u_n)u_n}{r+1}\right)\right|\\[3mm]
 
 &=&\left|\ds\int_{A_n} \frac{r}{r+1}\underline u^{r+1}+
 \ds\int_{B_n} \frac{r}{r+1}\underline u^{r}u_n\right|\\[3mm]
 
  &\leq& c_1+ c_2 \|u_n\|
 
 \end{array}
 \end{equation}
 and
   \begin{equation}\label{PSbdd3}
  \begin{array}{lll}
 \ds\left|\int_\Om G(x,u_n)-\frac{g(x,u_n)u_n}{r+1}\right|&=&\left|\ds\left(\int_{A_n} G(x,u_n)-\frac{g(x,u_n)u_n}{r+1}\right)+
 \ds\left(\int_{B_n} G(x,u_n)-\frac{g(x,u_n)u_n}{r+1}\right)\right|\\[3mm]
 
 &\leq &\left|\ds\int_{A_n} \left[\left(\frac{1}{q+1}-\frac{1}{r+1}\right) u_n^{q+1}-u_n+\frac{q}{q+1}\underline u^{q+1}+\frac{u_n}{r+1}\right] \right| \\[3mm] 
 
 && +\ds  \left| \ds\int_{B_n} \frac{r}{r+1}(\underline u^{q}-1)u_n \right|.\\[3mm]

 &\leq& c_3+c_4 \|u_n\|^{q+1}+ c_5 \|u_n\|
 \end{array}
 \end{equation}
 Now from (\ref{PSbdd1})-(\ref{PSbdd3}), clearly
 $$c+o_n(1) \geq  \left( \frac{1}{2}-\frac{1}{r+1}\right)\|u_n\|^2+C_1-C_2 \|u_n\|- C_3 \|u_n\|^{q+1}.$$
Since $q+1<2$ the above inequality implies that the Palais-Smale sequence is bounded in $X^s_0(\Omega).$ Thus $\{u_n\}$ admits a weakly convergent sub-sequence in $\xs.$ Since  $r<\ts,$ by a standard argument one can prove that the weakly convergent subsequence converges strongly in $\xs$ and hence the required result.  \hfill\qed. 
\begin{remark}
Using the same argument as above  one can  show that the Palais Smale sequence at any level $c$ is bounded even in case of $r=\tss.$
\end{remark}

\begin{theorem}{(Sub-Critical case)}
For $1<r<\ts-1$ and for each $\la>\la_0$ there exists a $\mu_\la>0$ such that the problem $(P_\la^\mu)$ admits at least two positive solutions whenever  $0<\mu<\mu_\la.$
\end{theorem}
\noi Proof: For a fixed $\la>\la_0,$ we prove the existence of first positive solution for $(P_\la^\mu)$ for all $\mu< \mu_{\la}$ via direct minimization method. Let $B_\rho$ denotes $\{u\in \xs : \|u\|<\rho\}.$ Using the estimate $\eqref{estFG}$ we find 
\begin{equation} \label{lower}
\ds I_\mu(u) \geq \frac{\|u\|^2}{2} -(\la+\mu)c|\Om| -(\mu+\la)c'\|u\|_{L^1}-\ds\la \|u\|_{{q+1}}^{q+1}-\mu \|u\|_{r+1}^{r+1}.
\end{equation}
 First we fix a $\rho$ large enough so that $\frac{\|u\|^2}{2} -\la c|\Om| -\la c'\|u\|_{L^1}-\ds\la \|u\|_{{q+1}}^{q+1} > 0$ on $\pa B_{\rho}.$ Next we choose $\mu$ small enough so that the right hand side of $\eqref{lower}$ is strictly positive on $\pa B_\rho.$
Thus for each $\la>\la_0,$ there exists a $\rho>0$ and $\mu_{\la}$ small such that 
\begin{equation}
\inf_{u\in \pa B_\rho} I_{\mu}(u) > 0 \mbox{ for all  } 0<\mu <\mu_\la. 
\end{equation} 
 Let $m_0=\inf_{u\in B_{\rho}} I_\mu(u).$  Since $I_\mu(0)=0$ we have $m_0\leq 0.$ 
It can easily be checked that $I_\mu$ is weakly lower semicontinuous in $B_\rho$ and hence attains a minimizer $u_0$ such that $I_\mu(u_0)=m_0.$ 
{The local minimizer $u_0$ is a weak solution of $(\tilde{P}_{\la}^\mu)$ and by Remark \ref{rem4.4} it is indeed a positive weak solution for $({P}_{\la}^\mu).$}

{
Now note that $I_\mu(t \phi_1) \rar -\infty$ as $t\rar \infty$ and due to the subcritical nature of $F,$ the functional $I_\mu$ satisfies the Palais- Smale condition at every level $c\in \R.$ Now choose a $t_0>>1$ so that $t_0 \phi_1 \not \in B_\rho$ and $I_\mu(t_0 \phi_1) < 0.$ Define the mountain pass critical level 
$$m_\mu = \ds \inf_{\gamma\in \Gamma} \, \max_{t\in [0,1]} I_\mu(\gamma(t)) $$
 where $\Gamma =\{ \gamma\in C([0,1]; \xs) : \ga(0)=0 \mbox{ and } \ga(1)= t_0 \phi_1 \}.$ For $0<\mu<\mu_\la$ clearly $m_\mu>0$ and hence $(P_\la^\mu)$ admits a second solution given by $\tilde{u}$ where $I_\mu(\tilde{u})=m_\mu.$}  \hfill\qed
 
Now we will prove the existence of one positive solution for $(P_\la^\mu)$ when $r=\tss.$

\begin{theorem}{(Critical case)}
For $r=\ts-1$ and for each $\la>\la_0$ there exists a $\mu_\la>0$ such that the problem $(P_\la^\mu)$ admits at least one positive solution whenever  $0<\mu<\mu_\la.$
\end{theorem}
\noi Proof:  As in the sub-critical case we obtain the first solution via the standard minimization technique. Now due to the presence of critical exponent term, the functional $I_\mu$ is no more weakly lower semicontinuous.  Let $\rho$ and $\mu_\la$ be as in $\eqref{lower}$ and  $\tilde{m}=\inf_{u\in B_{\rho}} I_\mu(u).$  Since $I_\mu(0)=0$ we have $\tilde{m}\leq 0.$ Now let $\{u_n\}$ denote a minimizing sequence of $I_\mu$ in $B_\rho$ such that 
\begin{equation}\label{unprop}
\left\{ \begin{array}{rll}
u_n   \rightharpoonup u_0 & \mbox{ weakly in } \xs,\\
 u_n   \rightharpoonup u_0 & \mbox{ and  weakly in } L^{\ts}(\Om) ,\\
 u_n   \rar u_0 & \mbox{ a.e  in }\Om \\
 \mbox{ and also } & \|u_n^+\|_{L^{\ts}} \mbox{ converges } . 
\end{array}\right.
\end{equation}
We write $I_\mu(u)= \ds L_\mu(u)- \la \int_{\Om} G(x,u)$ where 
\begin{equation}
L_\mu(u) = \frac{1}{2}\|u\|^2-\mu\int_{\Om} F(x,u).
\end{equation}
Using the strong convergence of $u_n$ in $L^{q+1}(\Om)$ it readily follows that $\ds \int_{\Om} G(x,u_n) \rar \int_{\Om} G(x,u_0)$ as $n\rar \infty.$ Now if we prove
 $L_{\mu}(u_0)\leq \liminf_{n\rar \infty} L_\mu(u_n),$ then $u_0$ is a minimizer of $I_\mu$ and hence a weak solution of $(P_\la^\mu).$

Thus it remains to prove that  for $\mu<<1,$ we have $\ds \liminf_{n\rar \infty} \left( L_{\mu}(u_n)-L_\mu(u_0)\right) \geq  0.$ 
The proof uses  idea  from \cite{MB} and \cite{Farkas} adapted for our functional $I_\mu$ which is cut off from below by $\ur$,
the details of which are given in the next Proposition \ref{prop1}. \hfill\qed. 

\begin{proposition}\label{prop1}
Assume that the sequence $\{u_n\}$ satisfies (\ref{unprop}). There exists a $\mu_{\la}>0$ such that  $\ds \liminf_{n\rar \infty} \left( L_{\mu}(u_n)-L_\mu(u_0)\right) \geq  0$ whenever $0<\mu<\mu_{\la}.$
\end{proposition}
\noi Proof: Given that 
\begin{equation}
L_\mu(u) = \frac{1}{2}\|u\|^2-\mu\int_{\Om} F(x,u).
\end{equation}
and $\{u_n\}$ is a minimizing sequence of $I_\mu$ in some fixed ball $B_R$ around origin and $\{u_n\}$ satisfies $\eqref{unprop}.$ We split the domain $\Om$ into four parts and write 
\begin{eqnarray}\label{Fidentity}
\ds \int_{\Om} F(x,u_0)-F(x,u_n) &= & \ds \left(\int_{A_1}+\int_{A_2}+\int_{A_3}+\int_{A_4}\right) \left( F(x,u_0)-F(x,u_n)\right)
\end{eqnarray}
for $A_1=\{x\in \Om: u_0\leq \ur, u_n\leq \ur\}, A_2=\{x\in \Om: u_0>\ur, u_n\leq\ur\},$  $A_3=\{x\in \Om: u_0\leq \ur, u_n>\ur\}$ and $A_4=\{x\in \Om: u_0>\ur, u_n>\ur\}.$ Note that each of the $A_i$ depends on $n.$ By dominated convergence theorem 
\begin{equation}
\left\{ \begin{array}{rll}
\ds \lim_{n\rar \infty} \ds \int_{A_i} \left( F(x,u_0)-F(x,u_n)\right) & = 
&0 \mbox{ for } i=1,2.\\
\ds \int_{A_3}\left( F(x,u_0)-F(x,u_n)\right)& = & \ds\int_{A_3}  \ds ( F(x,u_0)- \frac{u_n^{\ts}}{\ts}-\frac{\tss}{\ts}\ur^{\ts})\\[4mm]
\ds \int_{A_4}\left( F(x,u_0)-F(x,u_n)\right) &=  &\ds\frac{1}{\ts}\int_{A_4}(u_0^{\ts} -u_n^{\ts})
\end{array}\right.
\end{equation} 
Substituting above in $\eqref{Fidentity}$ we have
$$\begin{array}{lll}
\ds \liminf_{n\rar \infty} \ds \int_{\Om} F(x,u_0)-F(x,u_n)\\
\,\,\,\,\,\, \;\;\;\;\;\;\;\;\;\;\;\; \;\;\;\;\;\;\;\;=\ds \liminf_{n\rar \infty}
 \ds \left\{  \int_{A_3}  \ds ( F(x, u_0)- \frac{(u_0^+)^{\ts}}{\ts}-\frac{\tss}{\ts}\ur^{\ts})+
  \ds\int_{A_3\cup A_4}\frac{ (u_0^+)^{\ts}-u_n^{\ts}}{\ts} \right\}
\end{array}$$
\begin{equation} \label{eqn212a}
\end{equation}
Now note that $\chi_{A_3}$ converges to the characteristic function of $\{x: \ur=u_0\}.$ Thus $$\ds \lim_{n\rar \infty}\ds\int_{A_3}  \ds ( F(x,u_0)- \frac{u_0^{+ \ts}}{\ts} -\frac{\tss}{\ts}\ur^{\ts}) =0. $$ Next by dominated convergence theorem observe that $\ds \int_{A_1\cup A_2}(u_0^{+ \, \ts}-u_n^{+\, \ts}) \rar 0$ as $n\rar \infty.$ Thus from $\eqref{eqn212a}$ 
\begin{eqnarray}
\ds \liminf_{n\rar \infty} \ds \int_{\Om} F(x,u_0)-F(x,u_n)&=& \ds \liminf_{n\rar \infty} \frac{1}{\ts}\int_{\Om} (u_0^{+ \ts}-u_n^{+ \ts}) 
\end{eqnarray}
By Brezis-Lieb lemma $ \ds \lim_{n\rar \infty}\int_{\Om} (u_0^{+ \ts}-u_n^{+ \ts})$ exists and is equal to $ \ds -\lim_{n\rar \infty} \int_{\Om} |u_0^+-u_n^+|^{\ts}.$ Thus 
\begin{equation}
 \liminf_{n\rar \infty} \ds \int_{\Om} F(x,u_0)-F(x,u_n)=- \lim_{n\rar \infty} \frac{1}{\ts} \int_{\Om} |u_0^+-u_n^+|^{\ts}.
\end{equation}
Now we complete the proof of our claim following the proof of \cite{MB}. 
$$ \begin{array}{lll}\ds \liminf_{n\rar \infty} \left( L_{\mu}(u_n)-L_\mu(u_0)\right) &=&  \ds \liminf_{n\rar \infty} \left\{ \frac{\|u_n\|^2-\|u_0\|^2}{2}- \mu\int_{\Om}( F(x,u_n)-F(x,u_0)) \right\}\\
&\geq &  \ds \liminf_{n\rar \infty}  \frac{\|u_n\|^2-\|u_0\|^2}{2}+ \mu \liminf_{n\rar \infty} \ds \int_{\Om} F(x,u_0)-F(x,u_n) \\
&=& \ds \liminf_{n\rar \infty}  \frac{\|u_n\|^2-\|u_0\|^2}{2}- \frac{\mu}{\ts} \lim_{n\rar \infty}  \int_{\Om} |u_0^+-u_n^+|^{\ts}\\
&=&\ds \liminf_{n\rar \infty}\left\{ \frac{\|u_n\|^2-\|u_0\|^2}{2}- \frac{\mu}{\ts}  \int_{\Om} |u_0^+-u_n^+|^{\ts}\right\}\\
&\geq& \ds \liminf_{n\rar \infty} \left\{ \frac{1}{4}\|u_n-u_0\|^2-\frac{\mu}{\ts}\|u_n-u_0\|^{\ts}_{L^{\ts}}\right\}.
\end{array}$$
In the last step we have used equation $(3.8)$ of \cite{MB} which easily follows from the fact that $u_n \rightharpoonup  u_0$ and the inequality $b^2-a^2 \geq 2a(b-a)+\frac{1}{2}(a-b)^2.$
By Sobolev embedding $\ds \|u_n-u_0\|^{\ts}_{L^{\ts}} \leq S^{\ts}\|u_n-u_0\|^{\ts}$ and the  minimizing sequence $\{ u_n \}$ is bounded by $\rho$ independent of $\mu.$ 
$$ \begin{array}{lll}\ds \liminf_{n\rar \infty} \left( L_{\mu}(u_n)-L_\mu(u_0)\right) &\geq &
\ds \liminf_{n\rar \infty} \|u_n-u_0\|^2\left\{ \frac{1}{4} -\frac{\mu}{\ts} S^{\ts} \|u_n-u_0\|^{\ts-2} \right\}\\
&\geq &  \ds \liminf_{n\rar \infty} \|u_n-u_0\|^2 \left\{ \frac{1}{4} -\frac{\mu}{\ts} S^{\ts} (2\rho)^{\ts-2} \right\}.
\end{array}$$
Now the proposition follows if we choose $\ds \mu_0 \leq \frac{\ts}{4 S^{\ts} (2\rho)^{\ts-2}}.$\hfill\qed\\
\textbf{Proof of the main Theorem 1.2 :}  follows from theorems 4.4 and 4.5. 
\begin{remark}
 In case of a semilinear problem involving critical exponents, the assumption which is usually made to prove the compactness of the Palais Smale sequence is that of 
Ambrosetti-Rabinowitz in ({\cite{AR}}). More precisely if we define $H(t)=\int_0^th(s)ds$, (AR) condition reads as 
\begin{center}{{ there exists a constant } $a> 0$ {  such that
for }$|t|> a, \; 0< \theta H(t) < h(t)t$ where $\theta>2.$}\
\end{center}
Here (AR) condition is not satisfied for large negative values for $(\tilde{P}_\lambda^\mu)$ and hence we could not obtain a second solution for the critical exponent problem. 
\end{remark}
\noi \textbf{ACKNOWLEDGMENT:} Author Dhanya R  wishes to thank IIT Guwahati for their warm hospitality. She was supported by INSPIRE faculty fellowship  (DST/INSPIRE/04/2015/003221) when the work was being carried out.

\end{document}